\documentclass{amsart} 
\usepackage{amsmath,amscd} 
\usepackage{latexsym,amsmath,amssymb,mathrsfs,setspace,enumerate}
\usepackage[T1]{fontenc}
\usepackage[utf8]{inputenc}
\usepackage{lmodern} 
\usepackage{amssymb}
\usepackage{enumitem}
 \usepackage{tikz-cd}
\usepackage{setspace}
\usepackage{graphics}
\usepackage{graphicx}
\usepackage[utf8]{inputenc}
\usepackage[english]{babel}
\usepackage{fancyhdr}
\usepackage{tikz}
\usetikzlibrary{decorations.markings,positioning}
\theoremstyle{plain}

\newtheorem{lemma}{Lemma} 
 
\newtheorem{thm}{Theorem}
\newtheorem{exam}{Example}

\newtheorem{Rem}{Remark}
\newtheorem{defn}{Definition}

\tikzcdset{%
	triple line/.code={\tikzset{%
			double distance = 3pt, 
			double=\pgfkeysvalueof{/tikz/commutative diagrams/background color}}},
	quadruple line/.code={\tikzset{%
			double distance = 5.3pt, 
			double=\pgfkeysvalueof{/tikz/commutative diagrams/background color}}},
	Rrightarrow/.code={\tikzcdset{triple line}\pgfsetarrows{tikzcd implies cap-tikzcd implies}},
	RRightarrow/.code={\tikzcdset{quadruple line}\pgfsetarrows{tikzcd implies cap-tikzcd implies}}
}

\title{Ideal lattices of semigroup of doubly stochastic matrices}
\author {P. G. Romeo and Riya Jose}
\address{Dept. of Mathematics, Cochin University of Science and Technology, Kochi, Kerala, INDIA.}
\email{$romeo_-parackal@yahoo.com,\, riyajosemarangattu@gmail.com $}
\thanks{} 
\date{}

\begin{document}

\begin{abstract}
	In this paper we illustrate the rule for finding number of idempotents in the doubly stochastic matrix $D_n$ and also locate the idempotents in the case of $D_3$ and $D_4$. Further describe idempotent generated ideals of these semigroups and it is shown that the idempotent generated ideals of $D_n$ is a lattice. 
\end{abstract}
\maketitle

A real $n \times n$ matrix $A = [a_{ij}]$ is stochastic if it is non negative and  $\Sigma_{j=1}^n a_{ij} = 1$ for all $i$. The matrix $A = [a_{ij}]$ for which and both $\Sigma_{i=1}^n a_{ij} = 1$ for all $j$ and $\Sigma_{j=1}^n a_{ij} = 1$ for all $i$, then it is doubly stochastic. The fact that stochastic matrices readily arise in the study of homogenous Markov chains with finite number of states accounts for a large portion of the interest in this field of study. An $n$ state homogeneous Markov chain's state transition matrix is, in fact, an $n \times n$ stochastic matrix, and vice versa. If $A$ is an $n \times n$ stochastic matrix, there is an $n$ state homogeneous Markov chain whose state transition matrix is $A$. It is evident that the set of all $n \times n$ stochastic matrices $St_n$ is a semigroup under matrix multiplication and also the $n \times n$ doubly stochastic matrices $D_n$ is a sub semigroup of $St_n$. In the following we discuss the semigroup properties of $D_n$ and also investigate the structure of idempotents  generated ideals of $D_n$.

\section{Preliminaries}
Here we recall some definitios and results in semigroup theory and lattice theory. A semigroup $(S,\circ)$ is a non empty set $S$ together with an associative binary operation $\circ$ on it. The set of positive integers with multiplication as the operation is an example of a semigroup. 

\begin{defn}
	An element $s$ in a semigroup $( S,\circ)$ such taht If $s\circ s = s$ is referred to as an idempotent element of the semigroup $(S, \circ )$
\end{defn}
In the following for elements $x.y$ in a semigroup we simply write $xy$ for $x\circ y$ whenever there is no ambiguity regarding composition. A subset $I$ of $S$ that contains $IS =\{is | i \in I \,\text{and}\, s \in S\}$ or $SI =\{si | i \in I \,\text{and}\, s \in S\}$ is an ideal of a semigroup $S$. An ideal is said to be principal ideal if it is generated by a single element. 
The mathematical fields of order theory and abstract algebra both study lattices, which are abstract structures consisting of a partially ordered set in which each pair of elements has a distinct supremum (also known as a least upper limit or join) and a distinct infimum (also known as a maximum lower bound or meet).

Next we see a quick recap of several key concepts related to  semigroups of stochastic matrices and doubly stochastic matrices. A stochastic vector is a vector with non negative entries that add up to 1.It can be also defined using stochastic simplex as follows:
\begin{defn} (cf.\cite{PS}).
	A stochastic p- simplex $\Delta_p$ of $\mathbb{R}^n$ is a $(p-1) -$ simplex whose vertices are any p coordinate points of $\mathbb{R}^n$.
	
	$$\Delta_p = \{(a_1,a_2, ..., a_p): a_i \geq 0, \Sigma_{i=1}^p a_i = 1\}$$
and vectors belonging to $\Delta_p$ are called stochastic vectors.
\end{defn}
\begin{defn} A $ n \times m$ matrix $[a_{ij}]$ is said to be stochastic if all row vectors are of $[a_{ij}]$ are stochastic vectors.
\end{defn}
\begin{exam} The matrix 
	$	A = \begin{bmatrix}
		1/2 & 1/2 & 0 \\
		1 & 0 & 0\\
		1/4  & 1/4 & 1/2\\
	\end{bmatrix}$ is stochastic. Also it can be seen that the product of two stochastic matrices is again a stochastic matrix.
\end{exam}

\begin{defn} A $ n \times m$ matrix $[a_{ij}]$ is said to be doubly stochastic if all row vectors and column vectors are of $[a_{ij}]$ are stochastic vectors.
\end{defn}
\begin{exam}
The matrix 	$	\begin{bmatrix}
		1/2 & 1/2 & 0 \\
		0 & 0 & 1\\
		1/2 &  1/2 & 0\\
	\end{bmatrix}$
is dubly stochastic. It is easy to note that the product of two doubly stochastic matrices is again a doubly stochastic matrix.
\end{exam} 

 We write $St_n[D_n]$ to denote the set of all $n\times n$ stochastic [doubly stochastic] matrices. Note that if the transpose of matrix stochastic matrix  $A$ is also stochastic, then $A$ is doubly stochastic. 
 
 \begin{Rem}
 	The set $P_n$ of all $n \times n$ permutation matrices is a subgroup of $D_n$ which is isomorphic to $S_n$, the symmetric group on $n$ letters. 
 \end{Rem}
 \begin{lemma} (cf. \cite{Far}).
 	If a matrix and its inverse are both dubly stochastic then they belong to  to $P_n - $  the set of all $n \times n$ permutation matrices.  \end{lemma}
 
  Note that the idempotent $r\times r$ doubly-stochastic matrix,
 	$A=(a_{ij})$ has all its entries $a_{ij}$  equal to the number $\frac{1}{r}$. Also  for any doubly stochastic matrix $D$, $ AD = A = DA$.
 
 \begin{thm}
 	The set $St_n[D_n]$ of all $n\times n$ stochastic [doubly stochastic] matrices forms a compact  Hausdorff semigroup under matrix multiplication with the $n \times n$ identity matrix $I_n$ as identity element.
 \end{thm}
\begin{defn}
	An $n\times n$ square matrix $A$ is called reducible if there is a permutation matrix $P$ such that 
	$$PAP^{-1}=\begin{pmatrix}
		A_1&0\\
		B&A_2
	\end{pmatrix} $$
	where $A_1$ and $A_2$ are square matrices of orders $s$ and $n-s$ respectively.	
\end{defn}
The  Lemmas and Theorem bellow follows quickly. 
\begin{lemma}
A doubly-stochastic matrix is either irreducible or completely
reducible into irreducible doubly-stochastic matrices.
\end{lemma}
\begin{lemma}
	There exists a unique irreducible idempotent $r\times r$ doubly-stochastic matrix, namely the matrix 
	$A=(a_{ij})$  with all $a_{ij}$  equal to the number $\frac{1}{r}$. It is also the only idempotent doubly-stochastic matrix of rank $1$.
\end{lemma}
\begin{thm}
	Any idempotent $I\in D_n$ is of the form $PIP^T=U$  where $P$ is a permutation matrix 
	and $U$ is a matrix of the form
	$$
	U =\begin{pmatrix}
		{Q_1}&0&\cdots &0\\
		0&Q_2&\cdots &0\\
		\cdot&\cdot&\cdots\\
		0&0&\cdots &Q_s
	\end{pmatrix}
	$$
	Here $Q_i$ is $r_i\times r_i$ square matrix with all elements equal to $\frac{1}{r_i}$ and 
	$r_1+r_2+\cdots +r_s=n$. Conversely: Every matrix of this form is an idempotent in $D_n$ and it is of rank $s$.
\end{thm}
\section{Idempotents in $D_n$}

	Suppose $E$ is an idempotent in $D_n$ with rank $k$. By above theorem we know there exists an $n\times n$ permutation matrix $P$ such that 
	$$ PEP^T = E_1\oplus E_2\oplus\cdots \oplus E_k$$
	where each $E_i$ is $ n_i\times n_i$ matrix having each entry equal to $\frac{1}{n_i} \,: \, n_1+\cdots + n_k=n$ and $n_1\geq n_2\geq\cdots n_k\geq 1$.\\ 

Consider the semigroup $D_n$ and $\lambda$ be a partition of $n$ which has $\rho_{\alpha}$ parts equal to $\alpha$ ($1 \leq \alpha \leq n$), that is.,  $\lambda$ is of the form $n = \Sigma_{\alpha = 1}^n \rho_{\alpha} \alpha $. The number of idempotents in $D_n$ is given by $$ \Sigma \frac{n!}{(1!)^{\rho_1}\rho_1! (2!)^{\rho_2}\rho_2! \cdots (n!)^{\rho_n}\rho_n!}$$ where the sum extends over all partitions $n = \Sigma \rho_{\alpha} \alpha $ of $n$ (see cf.\cite{Far}).\\
	
Through the following example we illustrate the calculation of the number of idempotents of $D_n$ and list the idempotents of $D_n$ for $n= 3,4$.
\begin{exam}
	Consider $D_3$, the partitions of $3$ viz.,  $3,\,1+1+1,\, 2+1$. \\
	For $\lambda = 3, \rho_1 = 0, \rho_2 =0, \rho_3 = 1.$
	There is one idempotent of type $(3)$, namely the matrix
	$$
	E_1=\begin{pmatrix}
		\frac{1}{3}&\frac{1}{3}&\frac{1}{3}\\
		\frac{1}{3}&\frac{1}{3}&\frac{1}{3}\\
		\frac{1}{3}&\frac{1}{3}&\frac{1}{3}
	\end{pmatrix}
	$$ which is a rank $3$ idempotent.\\
For $\lambda = 2+1, \rho_1 = 1, \rho_2 =1, \rho_3 = 0.$	 There are three different idempotents of the form $(2+1)$ viz., 
	$$
	E_{2}^1=\begin{pmatrix}
		\frac{1}{2}&\frac{1}{2}&0\\
		\frac{1}{2}&\frac{1}{2}&0\\
		0&0&1
	\end{pmatrix},\,
	E_{2}^2=\begin{pmatrix}
		\frac{1}{2}&0&\frac{1}{2}\\
		0&1&0\\
		\frac{1}{2}&0&\frac{1}{2}
	\end{pmatrix},\,
	E_{2}^3=\begin{pmatrix}
		1&0&0\\
		0&\frac{1}{2}&\frac{1}{2}\\
		0&\frac{1}{2}&\frac{1}{2}
	\end{pmatrix}
	$$
	which are rank $2$ idempotents.\\
	For $\lambda = 1+1+1, \rho_1 = 3, \rho_2 = 0, \rho_3 = 0.$ and there is a unique idempotent of the form $(1+1+1)$
	$$
	E_3=\begin{pmatrix}
		1&0&0\\
		0&1&0\\
		0&0&1
	\end{pmatrix}
	$$ which is a rank $1$ idempotent. Hence all idempotents of $D_3$ are listed.

\end{exam}
\begin{exam}
	Consider $D_4$, the partitions of $4$ are $4,\,1+1+1+1,\, 2+2, 2+1+1, 3+1$. \\
	There is one idempotent of type $(4)$, namely the matrix
	$$
	E_1=\begin{pmatrix}
		\frac{1}{4}&\frac{1}{4}&\frac{1}{4} & \frac{1}{4}\\
			\frac{1}{4}&\frac{1}{4}&\frac{1}{4} & \frac{1}{4}\\
				\frac{1}{4}&\frac{1}{4}&\frac{1}{4} & \frac{1}{4}\\\	\frac{1}{4}&\frac{1}{4}&\frac{1}{4} & \frac{1}{4}
	\end{pmatrix}
	$$ which is a rank $1$ idempotent.\\
	Rank $2$ idempotents in $D_4$ are of type $(2+2)$ and $(3+1)$. There are $3$ of type $(2+2)$ and $4$ of type $(3+1)$ viz.,\\
	$$
	E_2^1=\begin{pmatrix}
		\frac{1}{2}&\frac{1}{2}&0& 0\\
		\frac{1}{2}&\frac{1}{2}&0& 0\\
		0& 0&\frac{1}{2}&\frac{1}{2}\\	
			0& 0&\frac{1}{2}&\frac{1}{2}\\
	\end{pmatrix}, 
E_2^2=\begin{pmatrix}
	\frac{1}{2}&0&\frac{1}{2}& 0\\
	0&\frac{1}{2}&0& \frac{1}{2}\\
     \frac{1}{2}&0&\frac{1}{2}& 0\\
		0&\frac{1}{2}&0& \frac{1}{2}\\
\end{pmatrix},
E_2^3=\begin{pmatrix}
	\frac{1}{2}&0&0& \frac{1}{2}\\
	0&\frac{1}{2}& \frac{1}{2}&0\\
	0&\frac{1}{2}& \frac{1}{2}&0\\
\frac{1}{2}&0&0& \frac{1}{2}\\
\end{pmatrix};$$
$$
E_2^4=\begin{pmatrix}
	1&0&0& 0\\
	0&\frac{1}{3}& \frac{1}{3}&\frac{1}{3}\\
0&\frac{1}{3}& \frac{1}{3}&\frac{1}{3}\\
0&\frac{1}{3}& \frac{1}{3}&\frac{1}{3}\\
\end{pmatrix},
E_2^5=\begin{pmatrix}

	\frac{1}{3}& 0& \frac{1}{3}&\frac{1}{3}\\
		0&1&0& 0\\
	\frac{1}{3}& 0& \frac{1}{3}&\frac{1}{3}\\
	\frac{1}{3}& 0& \frac{1}{3}&\frac{1}{3}\\
\end{pmatrix},
E_2^6=\begin{pmatrix}
	
	\frac{1}{3}& \frac{1}{3}& 0& \frac{1}{3}\\
	\frac{1}{3}& \frac{1}{3}& 0& \frac{1}{3}\\
	0&0&1& 0\\
	\frac{1}{3}& \frac{1}{3}& 0& \frac{1}{3}\\
\end{pmatrix},
E_2^7=\begin{pmatrix}
	
	\frac{1}{3}& \frac{1}{3}& \frac{1}{3}&0 \\
	\frac{1}{3}& \frac{1}{3}& \frac{1}{3}&0 \\
	\frac{1}{3}& \frac{1}{3}& \frac{1}{3}&0 \\
	0&0&0& 1\\
	
\end{pmatrix}.
$$
There are $6$ idempotents of type $(2+1+1)$ that are of rank $3$, they are 
$$E_3^1=\begin{pmatrix}
	
	1& 0& 0&0 \\
	0&	1& 0& 0 \\
	0 &0& \frac{1}{2}& \frac{1}{2} \\
	0 &0& \frac{1}{2}& \frac{1}{2} \\
	
\end{pmatrix},
E_3^2=\begin{pmatrix}
	
	 \frac{1}{2}& \frac{1}{2}& 0&0 \\
\frac{1}{2}& \frac{1}{2}& 0&0 \\
	0 &0&1&0 \\
	0 &0&0&1\\
	
\end{pmatrix},
E_3^3=\begin{pmatrix}
	1&0&0&0\\
	0&\frac{1}{2}&0& \frac{1}{2}\\
	0&0&1&0 \\
0&\frac{1}{2}&0& \frac{1}{2}\\
	
\end{pmatrix},
E_3^4=\begin{pmatrix}
	1&0&0&0\\
	0&\frac{1}{2}& \frac{1}{2}&0\\
		0&\frac{1}{2}& \frac{1}{2}&0\\
0&0&0&1\\
	
\end{pmatrix},\\
$$
$$
E_3^5=\begin{pmatrix}
		\frac{1}{2}& 0&\frac{1}{2}&0\\
		0&1&0&0\\
	\frac{1}{2}& 0&\frac{1}{2}&0\\
	0&0&0&1\\
	
\end{pmatrix},
E_3^6=\begin{pmatrix}
	\frac{1}{2}& 0 &0 &\frac{1}{2}\\
	0&1&0&0\\
	0&0&1&0\\
	\frac{1}{2}& 0 &0 &\frac{1}{2}\\
\end{pmatrix}
$$
Rank $4$ idempotent is of type $(1+1+1+1)$, namely the matrix
$$
E_4=\begin{pmatrix}
	1&0&0&0\\
	0&1&0&0\\
	0&0&1&0\\
	0&0&0&1\\
\end{pmatrix}
$$
\end{exam}

\section{Idempotent generated Ideals of $D_n$}

Now we proceed to describe ideals generated by the idempotents in $D_3$ and $D_4$. For $D_3$, the only ideal generated by rank 1 idempotent $E_1$ is 
$E_1D_3 = \{ E_1\}$ and the ideals generated by rank 2 idempotents $E_2^1, E_2^2, E_2^3$ are :\\

$E_2^1 D_3 =\{ \begin{pmatrix}
	a&b&1-a-b\\
	a&b&1-a-b\\
	1-2a&1-2b& -1 + 2a +2b
\end{pmatrix}| a, b \in \mathbb{R} \quad and \quad 0 \leq a, b \leq 1\}$\\

$E_2^2 D_3 =\{ \begin{pmatrix}
	\frac{1-a}{2}&\frac{1-b}{2}&\frac{a+b}{2}\\
	a&b&1-a-b\\
	\frac{1-a}{2}&\frac{1-b}{2}&\frac{a+b}{2}
\end{pmatrix}| a, b \in \mathbb{R} \quad and \quad 0 \leq a, b \leq 1\}$\\
	
$E_2^3 D_3 =\{ \begin{pmatrix}
	a&b&1-a-b\\
	\frac{1-a}{2}&\frac{1-b}{2}&\frac{a+b}{2}\\
	\frac{1-a}{2}&\frac{1-b}{2}&\frac{a+b}{2}
\end{pmatrix}| a, b \in \mathbb{R} \quad and \quad 0 \leq a, b \leq 1\}$.\\
	
The ideal generated by rank 3 idempotent is $E_3D_3 = D_3$.	\\

In a similar way, for the semigroup $D_4$, the ideal generated by rank 1 idempotent $E_1$ is $E_1D_4 = \{E_1\}$. Ideals generated by rank 2 idempotents are :\\
$E_2^1 D_4 =\{ \begin{pmatrix}
	a&b&c&1-a-b-c\\
	a&b&c&1-a-b-c\\
	\frac{1}{2}-a &\frac{1}{2}-b &\frac{1}{2}-c & -\frac{1}{2}+a+b+c\\ 
\frac{1}{2}-a &\frac{1}{2}-b &\frac{1}{2}-c & -\frac{1}{2}+a+b+c\\
\end{pmatrix}| a, b, c \in \mathbb{R} \, and \, 0 \leq a, b, c \leq 1\}$\\
$E_2^2 D_4 =\{ \begin{pmatrix}
	a&b&c&1-a-b-c\\
	\frac{1}{2}-a &\frac{1}{2}-b &\frac{1}{2}-c & -\frac{1}{2}+a+b+c\\ 
		a&b&c&1-a-b-c\\
	\frac{1}{2}-a &\frac{1}{2}-b &\frac{1}{2}-c & -\frac{1}{2}+a+b+c\\
\end{pmatrix}| a, b, c \in \mathbb{R} \, and \, 0 \leq a, b, c \leq 1\}$\\
$E_2^3 D_4 =\{ \begin{pmatrix}
	\frac{1}{2}-a &\frac{1}{2}-b &\frac{1}{2}-c & -\frac{1}{2}+a+b+c\\
	a&b&c&1-a-b-c\\
	a&b&c&1-a-b-c\\
	\frac{1}{2}-a &\frac{1}{2}-b &\frac{1}{2}-c & -\frac{1}{2}+a+b+c\\
\end{pmatrix}| a, b, c \in \mathbb{R} \, and \, 0 \leq a, b, c \leq 1\}$\\
$E_2^4 D_4 =\{ \begin{pmatrix}
	a&b&c&1-a-b-c\\
	\frac{1-a}{3} &\frac{1-b}{3} &\frac{1-c}{3}&\frac{a+b+c}{3}\\
	\frac{1-a}{3} &\frac{1-b}{3} &\frac{1-c}{3}&\frac{a+b+c}{3}\\
	\frac{1-a}{3} &\frac{1-b}{3} &\frac{1-c}{3}&\frac{a+b+c}{3}\\
\end{pmatrix}| a, b, c \in \mathbb{R} \, and \, 0 \leq a, b, c \leq 1\}$\\
$E_2^5 D_4 =\{ \begin{pmatrix}
	\frac{1-a}{3} &\frac{1-b}{3} &\frac{1-c}{3}&\frac{a+b+c}{3}\\
	a&b&c&1-a-b-c\\
	\frac{1-a}{3} &\frac{1-b}{3} &\frac{1-c}{3}&\frac{a+b+c}{3}\\
	\frac{1-a}{3} &\frac{1-b}{3} &\frac{1-c}{3}&\frac{a+b+c}{3}\\
\end{pmatrix}| a, b, c \in \mathbb{R} \, and \, 0 \leq a, b, c \leq 1\}$\\
$E_2^6 D_4 =\{ \begin{pmatrix}
	\frac{1-a}{3} &\frac{1-b}{3} &\frac{1-c}{3}&\frac{a+b+c}{3}\\
	\frac{1-a}{3} &\frac{1-b}{3} &\frac{1-c}{3}&\frac{a+b+c}{3}\\
		a&b&c&1-a-b-c\\
	\frac{1-a}{3} &\frac{1-b}{3} &\frac{1-c}{3}&\frac{a+b+c}{3}\\
\end{pmatrix}| a, b, c \in \mathbb{R} \, and \, 0 \leq a, b, c \leq 1\}$\\
$E_2^7 D_4 =\{ \begin{pmatrix}
	\frac{1-a}{3} &\frac{1-b}{3} &\frac{1-c}{3}&\frac{a+b+c}{3}\\
	\frac{1-a}{3} &\frac{1-b}{3} &\frac{1-c}{3}&\frac{a+b+c}{3}\\
	\frac{1-a}{3} &\frac{1-b}{3} &\frac{1-c}{3}&\frac{a+b+c}{3}\\
	a&b&c&1-a-b-c\\
\end{pmatrix}| a, b, c \in \mathbb{R} \, and \, 0 \leq a, b, c \leq 1\}$\\

Ideals generated by rank $3$ idempotents are:\\

$E_3^1 D_4 =\{ \begin{pmatrix}
	a&b&c&1-a-b-c\\
	d&e&f&1-d-e-f\\
	\frac{1-a-d}{2} &\frac{1-b-e}{2} &\frac{1-c-f}{2}&\frac{-1+a+b+c+d+e+f}{2}\\
	\frac{1-a-d}{2} &\frac{1-b-e}{2} &\frac{1-c-f}{2}&\frac{-1+a+b+c+d+e+f}{2}\\
	\end{pmatrix}| a, b, c, d, e, f \in \mathbb{R} \, and \, 0 \leq a, b, c, d, e, f \leq 1\}$\\

$E_3^2 D_4 =\{ \begin{pmatrix}
	a&b&c&1-a-b-c\\
	a&b&c&1-a-b-c\\
	d&e&f&1-d-e-f\\
	1-(2a+d) & 1-(2b+e) & 1-(2c+f) & -2+(2a+2b+2c+d+e+f)\\
	\end{pmatrix}\\ | a, b, c, d, e, f \in \mathbb{R} \, and \, 0 \leq a, b, c, d, e, f \leq 1\}$\\

$E_3^3 D_4 =\{ \begin{pmatrix}
	a&b&c&1-a-b-c\\
		\frac{1-a-d}{2} &\frac{1-b-e}{2} &\frac{1-c-f}{2}&\frac{-1+a+b+c+d+e+f}{2}\\
	d&e&f&1-d-e-f\\
	\frac{1-a-d}{2} &\frac{1-b-e}{2} &\frac{1-c-f}{2}&\frac{-1+a+b+c+d+e+f}{2}\\
\end{pmatrix}| a, b, c, d, e, f \in \mathbb{R} \, and \, 0 \leq a, b, c, d, e, f \leq 1\}$\\

$E_3^4 D_4 =\{ \begin{pmatrix}
	a&b&c&1-a-b-c\\
	d&e&f&1-d-e-f\\
	d&e&f&1-d-e-f\\
	1-(a+2d) & 1-(b+2e) & 1-(c+2f) & -2+(a+b+c+2d+2e+2f)\\
\end{pmatrix} \\| a, b, c, d, e, f \in \mathbb{R} \, and \, 0 \leq a, b, c, d, e, f \leq 1\}$\\

$E_3^5 D_4 =\{ \begin{pmatrix}
	a&b&c&1-a-b-c\\
	d&e&f&1-d-e-f\\
	a&b&c&1-a-b-c\\
	1-(2a+d) & 1-(2b+e) & 1-(2c+f) & -2+(2a+2b+2c+d+e+f)\\
\end{pmatrix} | \\a, b, c, d, e, f \in \mathbb{R} \, and \, 0 \leq a, b, c, d, e, f \leq 1\}$\\

$E_3^6 D_4 =\{ \begin{pmatrix}
		\frac{1-a-d}{2} &\frac{1-b-e}{2} &\frac{1-c-f}{2}&\frac{-1+a+b+c+d+e+f}{2}\\
	a&b&c&1-a-b-c\\
	d&e&f&1-d-e-f\\
	\frac{1-a-d}{2} &\frac{1-b-e}{2} &\frac{1-c-f}{2}&\frac{-1+a+b+c+d+e+f}{2}\\
\end{pmatrix}| a, b, c, d, e, f \in \mathbb{R} \, and \\ 0 \leq a, b, c, d, e, f \leq 1\}$\\
and the ideal generated by rank 3 idempotent is $E_3D_3 = D_3$.	\\

It is easy to observe that idempotents of same rank are $\mathcal{D}$ related and the ideals generated by $k$ rank idempotents contains ideal generated by $k+1$ rank idempotents. Ideal genrated by rank $1$ idempotent is the trivial one and ideal generated by rank $n$ idempotent is $D_n$. \\

Next we proceed to describe meet and join of these ideals generated by idempotents in $D_n$ and thus estabishes that the ideals generated by idempotents in $D_n$ forms a lattice. Here the meet between two ideals is the intersection and the join of two idempotent generated ideals is the smallest idempotent generated ideal containing the two. The following example illustrate the way in which the meet and join are described. 

\begin{exam}
Consider $ E_3^1, E_3^2 \in E(D_4)$, then the meet and join of the ideals generated by these idempotents can be located easily by using type of identical rows of idempotents. Note that the idempotent $ E_3^1$ in $D_4$ and all matrices in the ideal generated by it have third and fourth row identical. So we denote the ideal $<E_3^1> $ as $I_{(3,4)}^3$. Similarly the idempotent $ E_3^2$ in $D_4$ and all matrices in the ideal generated by it have first and second row identical and so  $<E_3^2> $ is denoted by $I_{(1,2)}^3$. Thus matrices in the intersection of these two ideals have first two and last two rows identical. The idempotent with this property is $E_2^1$, and is denoted as $I_{(1,2),(3,4)}^2$ and thus $I_{(3,4)}^3 \wedge I_{(1,2)}^3 = I_{(1,2),(3,4)}^2 .$
Hence we can define 
$$<E_3^1> \wedge <E_3^2> = <E_2^1>$$
 The join of the two ideals generated by idempotents $E$ and $E'$, will be the ideal generated by the idempotent which have only identical rows which are common to both $E$ and $E'$. Since for $ E_3^1, E_3^2 \in E(D_4)$, there are no common identical rows, hence the join of $I_{(3,4)}^3$ and $I_{(1,2)}^3$ will be the ideal generated by the idempotent which have no identical rows. i.e., $E_4$. Thus we have 
 $$<E_3^1> \vee <E_3^2> = <E_4>.$$
Following are few more cases of join and meet of ideal generated by idempotents in $D_4$:
 $$<E_2^1> \wedge <E_2^2> = I_{(1,2),(3,4)}^2 \wedge I_{(1,3),(2,4)}^2 = I_{(1,2,3,4)}^1 = <E_1>.$$
 $$<E_2^1> \vee <E_2^2> = I_{(1,2),(3,4)}^2 \vee I_{(1,3),(2,4)}^2 = I_{(1)(2)(3)(4)}^4 = <E_4>.$$
 $$<E_2^5> \wedge <E_2^6> = I_{(1,3,4)}^2 \wedge I_{(1,2,4)}^2 = I_{(1,2,3,4)}^1 = <E_1>.$$
 $$<E_2^5> \vee <E_2^6> =  I_{(1,3,4)}^2  \vee I_{(1,2,4)}^2 = I_{(1,4)}^3 = <E_3^6>.$$
 $$<E_2^1> \wedge <E_3^1> = I_{(1,2),(3,4)}^2 \wedge I_{(3,4)}^3 = I_{(1,2),(3,4)}^2 = <E_2^1>.$$
 $$<E_2^1> \vee <E_3^1> =  I_{(1,2),(3,4)}^2  \vee I_{(3,4)}^3 = I_{(3,4)}^3 = <E_3^1>.$$
 
\end{exam}

\end{document}